\documentclass[11pt,reqno]{amsart}
\usepackage[T1]{fontenc}
\usepackage{lmodern,amssymb,microtype}
\usepackage[letterpaper,margin=1in]{geometry}
\usepackage[hidelinks,hypertexnames=false]{hyperref}
\numberwithin{equation}{section}
\newtheorem{theorem}{Theorem}[section]
\newtheorem{proposition}[theorem]{Proposition}
\newtheorem{corollary}[theorem]{Corollary}
\newtheorem{lemma}[theorem]{Lemma}
\theoremstyle{definition}

\newtheorem{example}[theorem]{Example}
\theoremstyle{remark}
\newtheorem{remark}[theorem]{Remark}
\newcommand{\C}{\mathbb C}
\newcommand{\R}{\mathbb R}
\newcommand{\Q}{\mathbb Q}
\newcommand{\Z}{\mathbb Z}
\newcommand{\spn}{\operatorname{span}}
\newcommand{\Tr}{\operatorname{Tr}}
\newcommand{\sq}{\mathcal C}
\newcommand{\II}{\mathrm{II}_1}
\title[Number-theoretic obstructions for commuting square subfactors]
  {Number-theoretic obstructions\\for commuting square subfactors}
\author{Dietmar Bisch}
\address[Dietmar Bisch]{Department of Mathematics, Vanderbilt University, Nashville}
\thanks{Dietmar Bisch was supported in part by U.S. Army Research Office grant W911NF-23-1-0026.}
\author{Remus Nicoar\u{a}}
\address[Remus Nicoar\u{a}]{Department of Mathematics, University of Tennessee, Knoxville}
\date{November 29, 2026}
\hypersetup{pdftitle={Number-theoretic obstructions for commuting square subfactors},pdfauthor={Dietmar Bisch and Remus Nicoară}}
\begin{document}
\raggedbottom

\begin{abstract}
We prove that if a nondegenerate commuting square with connected inclusion
graphs and compatible Markov traces produces a finite-depth
subfactor, then the four algebras in the commuting square admit matrix 
units in which all four inclusions have \emph{algebraic} coefficients.
Fixing the algebras in the quadrilateral, the inclusion matrices and 
an upper bound on the depth of the resulting subfactor, we show that 
there are only finitely many such commuting squares up to trace-preserving
isomorphism.
More importantly, one number field suffices to define these commuting
squares. For commuting square subfactors of arbitrary depth, we determine 
from the minimal polynomial of the Jones index the set of rational numbers 
that occur as traces of projections in the higher relative commutants. 
This yields a powerful obstruction to realizing a hyperfinite subfactor 
as a commuting square subfactor.

We apply our techniques to Hadamard (or spin model) subfactors and
show that every entry of a dephased complex Hadamard matrix is
algebraic if its associated Hadamard subfactor has finite depth. Thus, we
obtain a number-theoretic criterion for infinite depth. It suffices
that one entry of the dephased matrix be transcendental.
\end{abstract}
\maketitle

\section{Introduction}

The construction of interesting subfactors of the hyperfinite II$_1$ factor
is often based on solving a finite-dimensional linear algebra problem. For
instance, to produce irreducible hyperfinite subfactors with small Jones
indices, one well-known approach is the construction of quadrilaterals
of multi-matrix algebras satisfying the \emph{commuting square condition}
(see for example \cite{Schou}, \cite{GHJ}). This condition was introduced
by Popa in \cite{P83} and so-called \emph{spin model commuting squares}, 
associated to Hadamard matrices, were considered for the first time 
in \cite{Orth}. We recall the definition of a commuting square.

Let $A_{1,1}$ be a finite direct sum of matrix algebras over $\C$, equipped 
with a faithful normalized trace $\tau$. Let $A_{1,0},A_{0,1}$ be unital 
$*$-subalgebras of $A_{1,1}$, and let $A_{0,0}$ be the unital $*$-subalgebra 
$A_{1,0}\cap A_{0,1}$. On each subalgebra, we consider the restriction 
of $\tau$. The quadrilateral
\begin{equation}\label{eq:square}
 \sq=\left(\begin{array}{ccc}
 A_{1,0}&\subset&A_{1,1}\\
 \cup&&\cup\\[-2pt]
 A_{0,0}&\subset&A_{0,1}
 \end{array},\,\tau\right)
\end{equation}
is a {\it commuting square} if the trace-preserving conditional 
expectations satisfy $E_{A_{1,0}}E_{A_{0,1}}=E_{A_{0,0}}$. Equivalently, 
$A_{1,0}\ominus A_{0,0}$ and $A_{0,1}\ominus A_{0,0}$ are orthogonal in 
$L^2(A_{1,1},\tau)$, i.e., with respect to the inner product
$\langle x, y \rangle = \tau(y^*x)$, $x$, $y \in A_{1,1}$. 
We will assume the commuting square is {\it nondegenerate}, i.e.,
satisfies $\spn A_{1,0}A_{0,1}=A_{1,1}$ \cite{P}.

For a nondegenerate commuting square with compatible Markov traces and
connected inclusion graphs,
write $A_{i,j}$ for the algebra 
in row $i$ and column $j$ of the iterated Jones basic 
construction \cite{GHJ,J,JS,Schou}:
\[
\begin{array}{ccccccc}
\vdots &&\vdots &&\vdots\\
\cup &&\cup &&\cup\\
 A_{2,0}&\subset&A_{2,1}&\subset&A_{2,2}&\subset&\cdots\\
 \cup&&\cup&&\cup\\
 A_{1,0}&\subset&A_{1,1}&\subset&A_{1,2}&\subset&\cdots\\
 \cup&&\cup&&\cup\\
 A_{0,0}&\subset&A_{0,1}&\subset&A_{0,2}&\subset&\cdots .
\end{array}
\]
Each elementary quadrilateral in the above grid is a commuting square. 
The rows and columns are obtained by iterating the Jones basic construction:
\[
 A_{i,j+1}=\langle A_{i,j},e_j\rangle\quad(j\ge1),\qquad
 A_{i+1,j}=\langle A_{i,j},f_i\rangle\quad(i\ge1),
\]
where $e_j$ implements the conditional expectation onto $A_{i,j-1}$ and 
is common to all rows, and $f_i$ implements the expectation onto 
$A_{i-1,j}$ and is common to all columns. Put 
$A_{\infty,j}=(\bigcup_{i\ge0}A_{i,j})''$ in the trace representation. 
The von Neumann algebras $N=A_{\infty,0}$ and $M=A_{\infty,1}$ are 
hyperfinite $\II$ factors, and the Jones tower of
$N\subset M$ is given by $M_{j-1}=A_{\infty,j}$ for $j\ge0$. 
Similarly for the horizontal limit algebras. We call 
$A_{\infty,0} \subset A_{\infty,1}$ and
$A_{0, \infty} \subset A_{1, \infty}$ {\it commuting square subfactors}
\cite{BCinf}.
Recall that
the subfactor has \emph{finite depth} if one (and hence both) of its 
principal graphs is finite. We measure depth on the principal graph 
rooted at the trivial $N$--$N$ bimodule. Its depth $D$ is the greatest 
graph distance from this root to a vertex.

\begin{theorem}\label{th:main}
Fix four finite-dimensional $C^*$-algebras $A_{0,0}$, $A_{0,1}$, 
$A_{1,0}$, $A_{1,1}$, four inclusion matrices with compatible labels 
on their simple summands and an integer $D\ge1$. 
Consider the nondegenerate commuting squares
\eqref{eq:square} with these data and a faithful normalized trace whose
restrictions are compatible Markov traces on all four inclusions. 
Suppose the four inclusion graphs are connected and the column
II$_1$ factors $A_{\infty,0}$ and $A_{\infty,1}$ give a finite
depth subfactor $A_{\infty,0}\subset A_{\infty,1}$ whose principal
graph has depth at most $D$.

Then, there are finitely many such commuting squares up to
trace-preserving isomorphism. Moreover, there is a number
field $K\subset\C$, closed under complex conjugation, such that for 
each commuting square one can choose $*$-matrix units in all four 
corners of the quadrilateral for which every inclusion sends each 
matrix unit to a linear combination of matrix units with coefficients in $K$.
\end{theorem}

We sketch the proof of Theorem~\ref{th:main}. Put $R_j=N'\cap M_{j-1}$. 
The depth bound is equivalent to $1\in\spn R_De_DR_D$. By Ocneanu 
compactness, $R_j=A_{1,0}'\cap A_{0,j}$, so the depth condition amounts 
to a finite system of matrix equations in the first two rows of the grid. 
For fixed inclusion data this system is defined over the field of real
algebraic numbers, that is, the real roots of nonzero polynomials in
$\Q[t]$. Indeed, the normalized Markov trace weights solve a linear
system over the field generated by the Perron-Frobenius eigenvalue of an integer
matrix. In path matrix units, the coefficients of the inclusions and
Jones projections are obtained from these weights by arithmetic and
positive square roots.

We first prove that the initial commuting squares fall
into at most countably many isomorphism classes. The fixed data and depth 
bound allow only finitely many principal graphs, hence finitely many 
hyperfinite subfactors by Ocneanu rigidity \cite[Corollary~2.2]{Oc}
and Popa's finite-depth classification theorem \cite{P90}.
It is shown in \cite{K} that there are countably many isomorphism
classes of commuting squares producing a fixed hyperfinite finite
depth subfactor.
In fixed real coordinates, each isomorphism class is a compact orbit.
The lexicographically least points of these orbits form a semialgebraic
set containing one representative of each class. Each coordinate projection of this countable set is
a countable semialgebraic subset of $\R$. A semialgebraic subset of
the line is a finite union of points and intervals, so a countable one
is finite. The representative set is therefore contained in a finite
product of finite sets. Since it is defined over the real algebraic
numbers, its coordinates are algebraic and generate a common number field.

For a commuting square subfactor $N\subset M$ of index $\lambda$,
the Markov trace on the grid of quadrilaterals gives
$\lambda^m\tau(p)\in\mathcal O_{\Q(\lambda)}$, the ring of integers
of the number field $\Q(\lambda)$, for every projection
$p\in N'\cap M_{2m-1}$ (Lemma~\ref{lem:scaled-traces}).
Thus at integer index the traces of projections in the higher relative
commutants are rational. Note that in the finite-depth case, rationality 
was already used in \cite{Bisch94} to obtain obstructions 
for intermediate and reduced subfactors.

Theorem~\ref{th:rational} determines the set of rational numbers that 
are traces of projections in the higher relative commutants of a
commuting square subfactor $N \subset M$ of arbitrary depth. 
If $b$ is the gcd of the
nonleading coefficients of the minimal polynomial of $[M:N]$, the set
of these traces is $[0,1]\cap\Z[1/b]$. Every value in this set is already
the trace of a projection in a Temperley--Lieb algebra generated by the
Jones projections. This applies in particular to every finite-depth
subfactor, since Popa's canonical commuting square determines its 
standard invariant, hence also the normalized traces of
projections in the higher relative commutants and the Jones 
projections \cite{P90}.

Corollary~\ref{cor:ade-traces} gives the exact rational trace range 
at Jones indices below four. In particular, the only rational projection
traces for the $E_6$ and $E_8$ principal graphs are $0$ and $1$, since
neither Coxeter number ($12$ or $30$) is twice a prime power.

The simplest commuting squares arise from Hadamard matrices. These 
spin model commuting squares are obtained as follows \cite{Orth}. 
If $D_d$ is the diagonal 
algebra in $M_d(\C)$ and $U$ is unitary, $D_d$ and $UD_dU^*$ are 
orthogonal over the scalars precisely when $H=\sqrt d\,U$ is a 
complex Hadamard matrix, that is, $|H_{ia}|=1$ and $HH^*=dI$. 
For $d\ge2$, consider the spin model commuting square with
\[
 A_{0,0}=\C1,\qquad A_{1,0}=D_d,\qquad
 A_{0,1}=UD_dU^*,\qquad A_{1,1}=M_d(\C)
\]
and normalized matrix trace $\tau$ (see also \cite{BHJ,N}). We obtain 
hyperfinite subfactors 
$N_H=A_{\infty,0}\subset M_H=A_{\infty,1}$ as mentioned above,
called the {\it spin subfactor} or {\it Hadamard subfactor} associated
to $H$. A Hadamard matrix $H$ is \emph{dephased} if its first row 
and column consist of ones. Every Hadamard matrix is equivalent to
a dephased matrix (by permuting rows and columns and multiplying rows and
columns by scalars of modulus one).

For a dephased $H$, an ordering of the minimal projections in the two 
abelian corners of the commuting square lets us recover each entry 
of $H$ from the commuting square. 
With $p_i=e_{ii}$ and $q_j=Ue_{jj}U^*$ one has $H_{ij}=d^3\tau(p_iq_jp_1q_1)$. It follows from Theorem~\ref{th:main} that every entry of $H$ is algebraic if the spin subfactor has finite depth. Thus a single transcendental entry in a dephased matrix implies infinite depth (Corollary~\ref{cor:had}).

For each ordered pair of rows define $\xi_{ij}(a)=H_{ia}\overline{H_{ja}}$ and let $P_{ij}=d^{-1}\xi_{ij}\xi_{ij}^*$ be the projection onto its span. Two distinct pairs are adjacent in the \emph{profile graph} if their row quotients 
have nonzero inner product \cite[Definitions~2.11.13--2.11.14]{Jplan}.
For a component $C$, write $V_C$ for the span of its row quotients and $Q_C$ for the orthogonal projection onto $V_C$. In Jones' spin model these are the minimal projections of $N_H'\cap M_{H,1}$, where $M_{H,1}=\langle M_H,e_1\rangle$ \cite[Theorem~2.11.15]{Jplan}. Qualitative cyclotomicity of their entries follows at arbitrary depth from the duality and splitting-field results for association schemes \cite{JMN,Mu}. At finite depth the order of the concrete block matrix $W_H=(P_{ij})_{i,j=1}^d\in M_{d^2}(\C)$ specifies a common cyclotomic field. We write $\zeta_m=\exp(2\pi i/m)$.

\begin{theorem}\label{th:moment}
Let $H$ be a dephased complex Hadamard matrix of order $d$ whose spin subfactor has finite depth. Then $W_H$ has finite order $m_H$.
\begin{enumerate}
\item[\textup{(i)}] For every $r\ge1$, the ordinary trace $Z_r(H)=\Tr(W_H^r)$ is a sum of $d^2$ $m_H$-th roots of unity. Every Galois conjugate of $Z_r(H)$ has absolute value at most $d^2$.
\item[\textup{(ii)}] For every component $C$ and every $a,b$, $d(Q_C)_{ab}\in\Z[\zeta_{m_H}]$. In particular $\sum_{i:(i,k)\in C} H_{ia}\overline{H_{ka}}\in\Z[\zeta_{m_H}]$ for every $k,a$.
\end{enumerate}
\end{theorem}

As an application of Theorem~\ref{th:moment}(i), let $C$ be a real
symmetric conference matrix of order $d>2$. This means that
$C=C^{\mathsf T}$, its diagonal is zero, its off-diagonal entries are
$\pm1$, and $C^2=(d-1)I$. Then the complex Hadamard matrix $H=I+iC$
has an infinite-depth spin subfactor (Proposition~\ref{prop:conference}).
The proof uses only the conference equations and
a rational, nonintegral value of $Z_2(H)$.

Let $\mu_v$ be the vertex dimensions of the finite principal graph of the spin subfactor, normalized by $\mu_*=1$ at the root. Its global index is $G=\sum_{v\ {\rm even}}\mu_v^2$. The next result determines a cyclotomic field from 
the graph data alone.

\begin{corollary}[Global-index obstruction]\label{cor:global}
Under the hypotheses of Theorem~\ref{th:moment}, let $A_H$ be the Hopf image of the row-quotient representation and let $h$ be its exponent. Then
\[
 G=\dim_{\C} A_H,\qquad m_H\mid h\mid G^3 .
\]
In particular, every $Z_r(H)$ and every entry of every $dQ_C$ belongs to $\Z[\zeta_{G^3}]$. Any prime ramified in the field generated over $\Q$ by these quantities must divide $G$.
\end{corollary}

We recall the defining identity for the Hopf exponent $h$ in Section~\ref{sec:profile}. The second power trace can be computed directly from the entries of $H$:
\begin{equation}\label{eq:zintro}
 Z_2(H)=d^{-2}\sum_{i,k}\left|\sum_a(H_{ia}\overline{H_{ka}})^2\right|^2 .
\end{equation}
For related infinite-depth criteria using the angle operator and tensor 
powers of the Hadamard representation, see for instance \cite{Mont,MontJ}. 

The paper is organized as follows. In Section~\ref{sec:models} we prove 
Theorem~\ref{th:main} by writing the depth condition in the first two rows 
of the grid of quadrilaterals. In Section~\ref{sec:traces} we prove the 
rational trace formula and use it to obtain obstructions for a subfactor
to be realized as a commuting square subfactor.
In Section~\ref{sec:profile} we apply Theorem~\ref{th:main} to dephased Hadamard matrices, prove the cyclotomic restrictions on $W_H$ and $dQ_C$, and apply 
the second power trace to the order-four family and conference matrices.

\medskip
\noindent\textbf{AI Statement.}
OpenAI's ChatGPT assisted with the development and exposition of the
arguments and with checking references. The authors take responsibility
for the mathematical content of the paper.

\section{Finiteness and algebraicity of commuting squares}\label{sec:models}

\begin{proof}[Proof of Theorem~\ref{th:main}]
Write $B_j=A_{0,j}$ and $C_j=A_{1,j}$ for the first two rows of the grid. Ocneanu compactness \cite[Theorem~5.7.6]{JS}, \cite[Theorem~2.11]{BC} gives
\begin{equation}\label{eq:compactness}
 R_j=N'\cap M_{j-1}=C_0'\cap B_j,
\end{equation}
for $j\ge1$, where the commutant on the right is taken in $C_j$. In the path model of the principal graph, the ideal $\spn R_De_DR_D\subset R_{D+1}$ has one matrix block for each endpoint occurring at level $D-1$. The complementary blocks correspond precisely to the vertices of depth $D+1$ \cite[\S1.3]{P}, \cite[\S4]{JP}. Hence
\begin{equation}\label{eq:depth}
 \operatorname{depth}(N\subset M)\le D
 \quad\Longleftrightarrow\quad
 1\in\spn R_D e_D R_D .
\end{equation}
Put $L=\dim_{\C}B_{D+1}$. The ideal lies in $B_{D+1}$ and has a basis chosen among the products $xe_Dy$, $x,y\in R_D$. Thus in \eqref{eq:depth} it suffices to express $1$ as a sum of at most $L$ such products.

We now fix the two horizontal inclusions $B_0\subset B_1$ and $C_0\subset C_1$ in standard block form. Their traces are determined by the data. For a connected Markov inclusion with matrix $F$ (source summands indexing the rows), the traces $t_v$ of rank-one projections in the target blocks $M_{b_v}(\C)$ are determined by
\begin{equation}\label{eq:weights}
 F^{\mathsf T}F t=\|F\|^2t,\qquad \sum_v b_vt_v=1.
\end{equation}
The Perron-Frobenius eigenvalue is simple, so $t_v\in\Q(\|F\|^2)$. In path matrix units the inclusions extend paths by one edge, and the Jones projections have coefficients obtained from these weights by arithmetic and positive square roots \cite[Ch.~4]{GHJ}, \cite[Definition~3.17 and Proposition~3.19]{JP}. Thus both row towers through column $D+1$ have algebraic coefficients. Write $e_j^B$ and $e_j^C$ for the horizontal Jones projections in the two fixed models.

Let $v_j:B_j\to C_j$ be trace-preserving unital $*$-homomorphisms, $0\le j\le D+1$. In the fixed row towers the compatibility conditions are
\[
 v_{j+1}|_{B_j}=v_j\quad(0\le j\le D),\qquad
 v_{j+1}(e_j^B)=e_j^C\quad(1\le j\le D).
\]
The pair $(v_0,v_1)$ gives a nondegenerate commuting square precisely when, in addition,
\[
 E^{C_1}_{C_0}v_1=v_0 E^{B_1}_{B_0},\qquad
 \spn C_0v_1(B_1)=C_1 .
\]
Here $E^V_U:V\to U$ is the trace-preserving conditional expectation. The basic construction extends each such pair to the row towers \cite[\S2]{K}. The extension is unique since $B_{j+1}=\spn B_j e_j^B B_j$ and $v_{j+1}(x e_j^B y)=v_j(x)e_j^C v_j(y)$ for $x,y\in B_j$. The corner traces and the prescribed vertical multiplicities determine each of the two initial vertical inclusions up to trace-preserving isomorphism. Thus the Markov property of each vertical inclusion and the isomorphism class of its column limit depend only on the fixed data. We may assume that the traces are Markov on both vertical inclusions and that both column limits are $\II$ factors, since otherwise the family in Theorem~\ref{th:main} is empty.

Let $\mathcal X_D$ consist of the maps above with the prescribed initial vertical multiplicities, together with $x_\ell,y_\ell\in B_D$, $1\le\ell\le L$, satisfying
\[
 [v_D(x_\ell),C_0]=[v_D(y_\ell),C_0]=0,\qquad
                 \sum_{\ell=1}^{L}x_\ell e_D^B y_\ell=1 .
\]
By \eqref{eq:compactness} and \eqref{eq:depth}, the pairs $(v_0,v_1)$ occurring in $\mathcal X_D$ are exactly the commuting squares under consideration. In the fixed matrix units, the $*$-homomorphism, trace and commutator relations, the tower compatibility conditions and the displayed depth identity are polynomial equations in the real and imaginary parts of the coordinates. The expectations have algebraic coefficients, since they are orthogonal projections for the two trace inner products. The multiplicity of a source summand in a target block is the rank in that block of the image of a minimal diagonal matrix unit. The spanning condition is also a rank condition on the products of basis elements. It follows that $\mathcal X_D$ is semialgebraic over $\overline{\Q}\cap\R$.

We first show that the commuting squares under consideration fall into at most countably many isomorphism classes. By Ocneanu compactness, $\dim_\C R_j\le\dim_\C B_j$ for $1\le j\le D+1$, and the right-hand sides are fixed by the initial data. Thus the block sizes and inclusion multiplicities in this initial segment of the relative-commutant tower are uniformly bounded. Since the depth is at most $D$, this segment determines the principal graph, so only finitely many principal graphs can occur. For each graph there are finitely many standard invariants by Ocneanu rigidity \cite[Corollary~2.2]{Oc}. Since the column limits are 
hyperfinite and the subfactor has finite depth, 
it follows from Popa's finite-depth classification theorem \cite{P90}
that each such standard invariant determines the subfactor up to 
isomorphism. Thus there are only finitely many possible subfactors up to isomorphism. For each proper subfactor on this list, Kawahigashi's 
result \cite{K}, applied 
with rows and columns interchanged, gives at most countably many isomorphism classes of commuting squares producing it \cite[Corollary~3.6 and the following paragraph]{K}. At index one each connected horizontal inclusion matrix has norm one and hence is $(1)$ \cite[Corollary~2.8(b) and Remark~2.10]{BC}, with the two directions interchanged. Both horizontal inclusions are then isomorphisms, so either prescribed vertical inclusion determines the commuting square up to isomorphism. Labeling the simple summands splits an isomorphism class into at most finitely many classes.

We now deduce finiteness from countability. Let $Y=\pi(\mathcal X_D)$,
where $\pi$ keeps the two initial vertical embeddings. By
Tarski--Seidenberg, $Y$ is semialgebraic over $\overline{\Q}\cap\R$.
Two points of $Y$ are equivalent if label-preserving tracial
automorphisms of the two fixed horizontal inclusions intertwine their
vertical embeddings. This relation is semialgebraic over the same
field, since the automorphism and intertwining conditions are polynomial
in real coordinates.

Each equivalence class is an orbit of the product of the two
automorphism groups. These groups are compact because they are closed subgroups
of the unitary groups for the trace inner products on the ambient row
algebras. Each orbit is therefore compact and has a unique
lexicographically least point, obtained by successively minimizing its
coordinates. Let $S$ consist of those $y\in Y$ for which there is no
equivalent $z\in Y$ with $z<_{\mathrm{lex}}y$. Lexicographic order is
semialgebraic, so Tarski--Seidenberg shows that $S$ is semialgebraic
over $\overline{\Q}\cap\R$. It contains exactly one point per class.
Thus $S$ is countable and hence finite \cite[Ch.~2]{BCR}.

Finally, a finite semialgebraic set defined over the real algebraic numbers consists of points with real algebraic coordinates \cite[\S1.4]{BCR}. The coefficients of the representatives in $S$, those of the fixed horizontal inclusions, and their complex conjugates therefore generate a number field $K$. Transporting the matrix units by the trace-preserving isomorphisms proves the theorem.
\end{proof}

\begin{remark}[A bound on the largest corner]\label{rem:corner-bound}
The fixed-corner and fixed-inclusion-matrix hypotheses in Theorem~\ref{th:main} may be replaced by a bound $\dim_\C A_{1,1}\le n$, with all other hypotheses retained. Indeed, every corner then has dimension at most $n$, so only finitely many isomorphism types of corner algebras and finitely many unital inclusion matrices can occur. Applying Theorem~\ref{th:main} to these finitely many data sets gives finiteness of the commuting squares. Taking the compositum of the resulting number fields gives a common field depending only on $n$ and $D$.
\end{remark}

\section{Traces in the higher relative commutants}\label{sec:traces}

Let $N\subset M$ be a finite-index inclusion of $\II$ factors with Jones tower $M_{-1}=N\subset M_0=M\subset M_1\subset\cdots$. Write $\lambda=[M:N]$, $R_n=N'\cap M_{n-1}$, and let $\tau$ denote the normalized trace.

\begin{lemma}\label{lem:scaled-traces}
Suppose $N\subset M$ arises from a nondegenerate commuting square of
finite-dimensional algebras with connected inclusion graphs and
compatible Markov traces. Let $K=\Q(\lambda)$ and let $\mathcal O_K$
be its ring of integers. Then $\lambda^m\tau(p)\in\mathcal O_K$ for
every $m\ge1$ and every projection $p\in R_{2m}$. The traces of
projections in all higher relative commutants generate $K$ as a field.
\end{lemma}
\begin{proof}
Write $B_j=A_{0,j}$ and let $F$ be the inclusion matrix
of $B_0\subset B_1$, with source blocks indexing its rows. Let $t_i$ be
the trace of a minimal projection in block $i$ of $B_0$. The Markov
trace gives $FF^{\mathsf T}t=\lambda t$, and a minimal projection in
block $j$ of $B_{2m}$ has trace $\lambda^{-m}t_j$. Connectedness and
the trace normalization determine $t$ over $K$.

By Ocneanu compactness \eqref{eq:compactness}, $p\in B_{2m}$. As $p$
commutes with $N$, it commutes with each minimal projection $h_i$ in
block $i$ of $B_0\subset N$, and $E_N(p)=\tau(p)1$. Thus
$\tau(h_i p)=t_i\tau(p)$. Computing the same trace in $B_{2m}$ gives
\[
 \lambda^m\tau(p)t_i=\sum_j T_{ij}t_j,\qquad
 T_{ij}=\operatorname{rank}\big((h_i p)|_j\big),
\]
where we choose one $h_i$ per block. Since $t_i>0$, this identity
places $\lambda^m\tau(p)$ in $K$ and makes it an eigenvalue of the
integer matrix $T$. It therefore belongs to $\mathcal O_K$.

Every projection in $M_i'\cap M_j$ belongs to $R_{2m}$ for
$2m\ge j+1$, so its trace lies in $K$. Since the Jones projection
$e_1$ has trace $\lambda^{-1}$, these traces generate $K$.
\end{proof}

At finite depth the index field $\Q(\lambda)$ is a subfield of a
cyclotomic field \cite[Corollary~8.54]{ENO}, \cite[\S5.3.1]{JMS}.

\begin{theorem}[Rational projection traces]\label{th:rational}
Let $N\subset M$ be the subfactor arising from a nondegenerate commuting
square of finite-dimensional algebras with connected inclusion graphs and
compatible Markov traces. Write the minimal polynomial of $\lambda=[M:N]$ as
\[
 f(t)=t^s+c_{s-1}t^{s-1}+\cdots+c_0,\qquad b=\gcd(|c_0|,\ldots,|c_{s-1}|).
\]
The set of rational numbers realized as normalized traces of projections in the higher relative commutants $M_i'\cap M_j$, $-1\le i\le j$, is
\[
                      [0,1]\cap\Z[1/b].
\]
Here $\Z[1/1]=\Z$. Equivalently, a reduced fraction $a/q\in[0,1]$ occurs if and only if every prime divisor of $q$ divides $b$. Every such fraction is realized by a projection in the Temperley--Lieb algebra generated by $e_1,\ldots,e_{2m-1}$ in $N'\cap M_{2m-1}$ for some $m\ge1$.
\end{theorem}

\begin{proof}
Let $p$ be a projection in a higher relative commutant with
$\tau(p)=a/q$ in lowest terms. Choose $m\ge1$ such that $p\in R_{2m}$,
and write $\mathcal O$ for the ring of integers of $\Q(\lambda)$.
Lemma~\ref{lem:scaled-traces} gives
\[
 a\lambda^m\in q\mathcal O,\quad (a,q)=1
                 \quad\Longrightarrow\quad \lambda^m\in q\mathcal O .
\]
If $\ell$ is a prime divisor of $q$, multiplication by $\lambda$ on $\mathcal O/\ell\mathcal O$ is therefore nilpotent. Its characteristic polynomial is $f(t)\bmod\ell$, since multiplication by $\lambda$ on $\Q(\lambda)$ has characteristic polynomial $f$ and preserves $\mathcal O$. Thus $f(t)\equiv t^s\bmod\ell$, or equivalently $\ell\mid b$.

For the converse, the index-one case is immediate. Put $r=\lambda^{-1}$,
so $0<r\le1/2$ by the index theorem \cite{J}. We first describe a
collection of orthogonal projections using only the Jones projections.
The relations $e_1e_2e_1=re_1$ and $e_2e_1e_2=re_2$ give
\[
 v=\frac{(1-e_1)e_2e_1}{\sqrt{r(1-r)}},\qquad
 v^*v=e_1,\qquad vv^*\le 1-e_1.
\]
Thus $e_1,vv^*,1-e_1-vv^*$ are orthogonal projections of traces
$r,r,1-2r$. Form $n$ copies of this triple, using the pairs
$(e_{3h+1},e_{3h+2})$ for $0\le h<n$. Projections in different
triples commute by the Jones relations. Their traces multiply when we
choose one projection from each triple. Indeed, if $a$ is a product
from the preceding triples, the Markov property gives
$\tau(ae_{3h+1})=r\tau(a)$. Traciality and commutation with the
corresponding copy of $v$ give the same trace for the second choice,
and subtraction gives $(1-2r)\tau(a)$ for the third.

We claim that if $P\in\Z[x]$, $P(0)=0$, and $0<P(x)<1$ for
$0<x\le1/2$, then $P(r)$ is the trace of a sum of these mutually
orthogonal product projections for some $n$. The argument is a
one-variable partition-polynomial calculation (see also \cite{Yingst}).
For $0\le j\le n$, there are $C_{n,j}=2^j\binom nj$ choices in which
the first or second projection is used exactly $j$ times. Each has
trace $r^j(1-2r)^{n-j}$. Expand $P(x)$ in the polynomial basis
$x^j(1-2x)^{n-j}$, $0\le j\le n$. For the monomial $x^k$, expansion of
$x^k(2x+1-2x)^{n-k}$ gives the integer coefficient
$2^{j-k}\binom{n-k}{j-k}$ when $j\ge k$, and zero otherwise. Dividing
that coefficient by $C_{n,j}$ gives
\[
  2^{-k}\frac{(j)_k}{(n)_k}
  =\left(\frac{j}{2n}\right)^k+O_k(n^{-1})
\]
uniformly in $j$, where $(j)_k=j(j-1)\cdots(j-k+1)$ and $(j)_0=1$.
Consequently the normalized coefficients of any fixed polynomial $Q$
differ from $Q(j/(2n))$ by $O_Q(n^{-1})$. If $Q$ is strictly positive
on $[0,1/2]$, it follows that all its coefficients in this basis are
positive for $n$ large enough.

Write $P=x^d Q$ where $Q(0)\ne0$. Then $Q$ is strictly positive on
$[0,1/2]$. Expanding $Q$ at degree $n-d$ and multiplying by $x^d$
shows that the coefficients of $P$ at degree $n$ are nonnegative.
The same positivity argument applies directly to $1-P$. Since the
coefficients of $1=(2x+1-2x)^n$ are $C_{n,j}$, the $j$th coefficient
of $P$ is an integer between $0$ and $C_{n,j}$. Selecting that many
product projections of each type proves the claim.

It remains to construct such a polynomial for each allowed fraction
$0<a/q<1$. For each prime $\ell\mid q$ we have $\ell\mid b$, so
$f(\lambda)=0$ gives
\[
 \frac1\ell=-\sum_{j=0}^{s-1}\frac{c_j}{\ell}r^{s-j}.
\]
Multiplying these polynomials with the prime multiplicities of $q$
gives $U\in x\Z[x]$ with $U(r)=1/q$. Set
$D(x)=1-qU(x)$ and $H(x)=1-qx(1-x)^q$. On $[0,1/2]$ we have
$0<H(x)\le\exp(-q2^{-q}x)$. Indeed,
$qx(1-x)^{q-1}\le1$ by the binomial theorem and
$(1-x)^q\ge2^{-q}$. Since $D(0)=1$, some $C>0$ satisfies
$|D(x)|\le1+Cx\le e^{Cx}$ on this interval. Choose an integer $L$
with $Lq2^{-q}>C$ and put
\[
             P(x)=\frac aq\bigl(1-D(x)^2 H(x)^{2L}\bigr).
\]
Both $D$ and $H$ are congruent to $1$ modulo $q$, so $P\in\Z[x]$.
We have $P(0)=0$ and $P(r)=a/q$ since $D(r)=0$. The choice of $L$
ensures $|D(x)H(x)^L|<1$ for $0<x\le1/2$, whence
$0<P(x)\le a/q<1$. The claim gives a projection of trace $a/q$
generated by finitely many Jones projections. Taking $2m\ge3n$
places it in $R_{2m}$. The endpoints are the traces of $0$ and $1$.
\end{proof}

The realization by Temperley--Lieb projections in the converse uses
only the Jones relations, the Markov trace and the fact that $\lambda$
is an algebraic integer. Thus for any finite-index subfactor with
algebraic integer index $\lambda$, every number in $[0,1]\cap\Z[1/b]$
occurs as such a trace.

At integer index $d$, all projection traces in the higher relative
commutants of a commuting square subfactor are rational, since
$d^m\tau(p)\in\Z$ for $p\in R_{2m}$ by Lemma~\ref{lem:scaled-traces}. 
The idea of using rationality already appears in 
\cite[Theorem~3.2, p.~33]{Bisch94}, where it is shown
that for a finite-depth inclusion of separable $\II$ factors with
integer index, the Jones index of every intermediate and every reduced 
subfactor must also be an integer.
Theorem~\ref{th:rational}
determines exactly which traces occur. Since $f(t)=t-d$, their range
is $[0,1]\cap\Z[1/d]$.
At index six, for example, projections of traces $1/2$ and $1/3$ occur.

\begin{corollary}[Rational traces at indices below four]\label{cor:ade-traces}
Let $N\subset M$ be a subfactor with index
\[
 [M:N]=4\cos^2(\pi/n),\qquad n\in\Z,\quad n\ge3.
\]
The rational numbers realized as traces of projections in its higher relative commutants are exactly
\[
 \begin{cases}
 [0,1]\cap\Z[1/p],& n=2p^a\text{ for a prime }p\text{ and an integer }a\ge1,\\
 \{0,1\},&\text{otherwise}.
 \end{cases}
\]
In particular, every nonzero proper projection in a higher relative commutant of a subfactor with principal graph $E_6$ or $E_8$ has irrational trace.
\end{corollary}
\begin{proof}
Since $[M:N]<4$, the subfactor has finite depth \cite{JMS}.
Its standard invariant therefore has a canonical commuting square
realization \cite{P90}, to which Theorem~\ref{th:rational} applies.
Put $\lambda=4\cos^2(\pi/n)$, $K=\Q(\zeta_n)$ and $F=\Q(\lambda)=K^+$. Since $\lambda=(1+\zeta_n)(1+\zeta_n^{-1})$, norm transitivity and the standard evaluation of cyclotomic polynomials at $-1$ give
\[
 N_{F/\Q}(\lambda)=N_{K/\Q}(1+\zeta_n)=\Phi_n(-1)
 =\begin{cases}p,&n=2p^a,\\1,&\text{otherwise},\end{cases}
\]
where $p$ is prime and $a\ge1$ \cite{BHM}. Let $f$ and $b$ be as in Theorem~\ref{th:rational}. If this norm is one, the constant coefficient of $f$ is $\pm1$, so $b=1$. If $n=2p^a$, the prime $p$ is totally ramified in $K$, and $\zeta_n$ reduces to $-1$ at the unique prime above $p$ \cite[Ch.~2]{Wash}. Every conjugate of $\lambda$ therefore reduces to zero, so $f(t)\equiv t^{\deg f}\pmod p$. The constant coefficient of $f$ has absolute value $p$, hence $b=p$. Theorem~\ref{th:rational} proves the formula. The $E_6$ and $E_8$ graphs have Coxeter numbers $12$ and $30$, respectively, and both fall in the second case.
\end{proof}

\begin{example}[Rational traces for two $A_\infty$ commuting square subfactors]\label{ex:infinite-traces}
At the Haagerup and Asaeda--Haagerup indices $(5+\sqrt{13})/2$ and
$(5+\sqrt{17})/2$ \cite{AH}, the minimal polynomials are $t^2-5t+3$
and $t^2-5t+2$, respectively. Both have $b=1$. Thus the only rational
projection traces are $0$ and $1$ for every commuting square subfactor
at either index.

\medskip

Bisch and C\'aceres constructed a commuting square subfactor with
principal graph $A_\infty$ at the Asaeda--Haagerup index \cite{BCinf}.
The preceding calculation therefore excludes a projection of trace
$1/2$ in any of its higher relative commutants, although $2$ divides
the norm of the index.

For their commuting square subfactor with the same principal graph at
$\lambda=3+\sqrt3$ \cite{BCinf}, the minimal polynomial is $t^2-6t+6$.
Here $b=6$, so Theorem~\ref{th:rational} gives the rational projection
traces $[0,1]\cap\Z[1/6]$. To see $1/2$ explicitly, consider the four-strand
Temperley--Lieb algebra in $N'\cap M_3$. There are three paths of
length four from the root to the vertex at depth two, whose weight
is $\lambda-1$. The three diagonal minimal projections of the
corresponding $M_3(\C)$ block each have trace
$(\lambda-1)/\lambda^2=1/6$. The identity of that block has trace $1/2$.
\end{example}

\begin{example}[An obstruction to a commuting square 
presentation]\label{ex:presentation}
Let $A\subset B$ be the irreducible hyperfinite subfactor of index
$9/2$ constructed in \cite{Bisch94}, 
and $S\subset T$ an outer cyclic group inclusion
of index two. The tensor product
$N=A\mathbin{\bar\otimes} S\subset M=B\mathbin{\bar\otimes}T$
is irreducible and hyperfinite of index nine. The intermediate
$P=B\mathbin{\bar\otimes}S$ has $[M:P]=2$, so the projection
$e_P\in N'\cap M_1$ onto $L^2(P)$ has normalized trace $1/2$.
At index nine Theorem~\ref{th:rational} allows only denominators that
are powers of three. Consequently $N\subset M$ admits no presentation
by a nondegenerate finite-dimensional commuting square with connected
inclusion graphs and compatible Markov traces.

The same projection gives an obstruction at the irrational index
$\beta=9\alpha$, where $\alpha=(5+\sqrt{17})/2$. Tensor $N\subset M$
with the Asaeda--Haagerup inclusion $U\subset V$ of index $\alpha$
\cite{AH}. The intermediate $P\mathbin{\bar\otimes}V$ still has
upper index two, with Jones projection $e_P\otimes 1$ of trace $1/2$.
The minimal polynomial of $\beta$ is $t^2-45t+162$, so $b=9$.
Although $2$ divides the norm $162$, Theorem~\ref{th:rational}
excludes $1/2$ and rules out a commuting square presentation.
At the same index $\beta$, tensoring $U\subset V$ with an outer cyclic
group inclusion of index nine gives an irreducible hyperfinite
finite-depth subfactor, which has a finite-dimensional commuting square
presentation \cite{P90}.
\end{example}

\section{Hadamard matrices and spin model commuting squares}\label{sec:profile}

\begin{corollary}\label{cor:had}
If the spin subfactor of a dephased complex Hadamard matrix $H$ has finite depth, every entry of $H$ is algebraic. In particular, a single transcendental entry implies infinite depth.
\end{corollary}
\begin{proof}
Write $U=H/\sqrt d$, let $u_j$ be its $j$th column and let $\varepsilon_i$ be the $i$th coordinate vector. For $p_i=e_{ii}$ and $q_j=Ue_{jj}U^*$ one has
\[
 p_iq_j=U_{ij}\,\varepsilon_i u_j^*,\qquad
 \spn\{p_iq_j:1\le i,j\le d\}=M_d(\C).
\]
The spanning assertion follows because every $U_{ij}$ is nonzero and the operators $\varepsilon_i u_j^*$ form a basis of $M_d(\C)$. Thus the spin model commuting square is nondegenerate. Its four inclusions are connected Markov inclusions for $\tau=d^{-1}\Tr$, and $[M_H:N_H]=d$ \cite{JS,N}.

In a finite-dimensional abelian $C^*$-algebra, each simple summand is one-dimensional and its unit is its only matrix unit. Thus, up to ordering, the only systems of matrix units in $D_d$ and $UD_dU^*$ are $(p_i)_{i=1}^d$ and $(q_j)_{j=1}^d$, respectively. Theorem~\ref{th:main} therefore gives a single system of matrix units in $M_d(\C)$ in which every $p_i$ and every $q_j$ has coefficients in $K$. In these coordinates $p_iq_jp_1q_1$ belongs to $M_d(K)$, so its normalized trace lies in $K$. Hence
\[
 H_{ij}=H_{ij}\overline{H_{1j}}H_{11}\overline{H_{i1}}
       =d^3\tau(p_iq_jp_1q_1)\in K .
\]
The first equality uses dephasing. The middle expression is an element of the Haagerup set \cite{Haag}.
\end{proof}

\begin{proof}[Proof of Theorem~\ref{th:moment}]
The row quotients $d^{-1/2}\xi_{ij}$ form an orthonormal basis of $\C^d$ when either $i$ or $j$ is fixed. Thus $(P_{ij})$ is a magic unitary; in particular, $W_H$ is unitary. It defines a representation $\pi_H:u_{ij}\mapsto P_{ij}$ of the quantum permutation algebra. Write $A_H$ for its Hopf image, the smallest quantum group quotient through which $\pi_H$ factors \cite{BN,BBS}, and $u=(u_{ij})$ for its fundamental corepresentation. Then $W_H=(\mathrm{id}\otimes\pi_H)(u)$. The description of the standard invariant by the corresponding coaction 
can be found in \cite[Corollary~3.1]{BanK}.

Let $A$ be the adjacency matrix of the finite principal graph of $N_H\subset M_H$, rooted at $*$. Since $A_{1,0}=D_d$, Ocneanu compactness identifies $N_H'\cap A_{\infty,n}$ with $D_d'\cap A_{0,n}$ inside $A_{1,n}$. This is the $n$-box space $\mathcal B_n$ of Jones' row-quotient spin model \cite[\S2.11]{Jplan}, \cite[\S2]{N}. The fixed-space formula in \cite[Theorem~5.11(1)]{BanQ} and the 
path description of this relative commutant give
\[
 \dim\operatorname{Fix}(u^{\otimes n})=\dim\mathcal B_n=(A^{2n})_{**}.
\]
For $r\ge0$, let $v_r(V)$ be the multiplicity of an irreducible corepresentation $V$ in $u^{\otimes r}$. By self-duality of $u$,
\[
 \langle v_r,v_s\rangle
 =\dim\operatorname{Hom}(u^{\otimes r},u^{\otimes s})
 =(A^{2r+2s})_{**}=\langle A^{2r}\delta_*,A^{2s}\delta_*\rangle .
\]
The common Gram matrix has rank at most the number of vertices of the principal graph. Choose a basis for $\spn\{v_r:r\ge0\}$ among the $v_r$. The union of its supports is finite and contains the support of every $v_r$. Since $u$ is self-dual and its coefficients generate $A_H$, every irreducible occurs in some $u^{\otimes r}$. Thus there are finitely many irreducibles, and $\dim_\C A_H=\sum_V(\dim V)^2<\infty$ by the matrix-coalgebra decomposition.

The finite-dimensional Hopf $C^*$-algebra $A_H$ is semisimple and cosemisimple, with involutive antipode. Its exponent $h$ is the least positive integer such that $m^{(h)}\Delta^{(h)}(x)=\varepsilon(x)1$ for all $x\in A_H$, where $\Delta^{(h)}$ is the iterated coproduct with $h$ tensor factors and $m^{(h)}$ multiplies them in order. Etingof and Gelaki give $h\mid(\dim_\C A_H)^3$ \cite[Definition~2.1 and Theorem~4.3]{EG}. On the generators the exponent identity gives
\[
 (u^h)_{ij}=m^{(h)}\Delta^{(h)}(u_{ij})=\delta_{ij}1,
 \qquad W_H^h=I .
\]
Hence $m_H=\operatorname{ord}(W_H)$ divides $h$. Every $Z_r(H)$ is the sum of the $r$th powers of the $d^2$ eigenvalues of $W_H$. These are $m_H$-th roots of unity, as are their Galois conjugates. This proves (i).

For (ii), fix a profile component $C$ and put $B_C=(\mathbf 1_{(i,k)\in C})_{i,k}$. Distinct component spaces are orthogonal. Decomposing each row basis, or each column basis, gives
\begin{equation}\label{eq:component-sum}
 Q_C=\sum_{k:(i,k)\in C}P_{ik}\quad\text{for each }i,\qquad
 Q_C=\sum_{i:(i,k)\in C}P_{ik}\quad\text{for each }k.
\end{equation}
In particular $I\otimes Q_C$ commutes with $W_H$. The $(i,k)$ block of $W_H$ restricted to $\C^d\otimes V_C$ has trace $\mathbf 1_{(i,k)\in C}$. Its unnormalized partial trace over $V_C$ is therefore $B_C$.

Set $w_{ab}(i)=H_{ia}\overline{H_{ib}}$. Taking the $(b,a)$ entry of the first identity in \eqref{eq:component-sum} gives
\[
 (B_Cw_{ab})(i)=\sum_{k:(i,k)\in C}H_{ka}\overline{H_{kb}}
                  =d(Q_C)_{ba}\,w_{ab}(i).
\]
For fixed $b$ the $w_{ab}$ form an orthogonal basis, so $B_C$ is normal. Each spectral projection $E_\eta$ of $B_C$ commutes with $w_{ab}w_{ab}^*$ for every $a,b$. Writing $e_{ab}$ for the standard matrix units, we have
\[
 W_H=\sum_{i,k}e_{ik}\otimes P_{ik}
    =\frac1d\sum_{a,b} w_{ab}w_{ab}^*\otimes e_{ab}.
\]
The second expression gives $[E_\eta\otimes I,W_H]=0$. Since $I\otimes Q_C$ also commutes with $W_H$, the range of $E_\eta\otimes Q_C$ reduces $W_H$. The partial-trace identity gives
\[
 (\operatorname{rank}E_\eta)\eta
 =\Tr(E_\eta B_C)
 =\Tr\bigl(W_H|_{\operatorname{ran}E_\eta\otimes V_C}\bigr)
 \in\Z[\zeta_{m_H}].
\]
The restricted unitary has order dividing $m_H$, which gives the last inclusion. Dividing by $\operatorname{rank}E_\eta$ gives $\eta\in\Q(\zeta_{m_H})$. Since $B_C$ has integer entries, $\eta$ belongs to the ring of integers $\Z[\zeta_{m_H}]$ \cite[Ch.~2]{Wash}. Apply this to $\eta=d(Q_C)_{ba}$. The last assertion of (ii) follows from the $(a,1)$ entry of the second identity in \eqref{eq:component-sum} and the dephasing of $H$.
\end{proof}

\begin{proof}[Proof of Corollary~\ref{cor:global}]
For irreducible corepresentations $V,W$ of $A_H$, put $F_{VW}=\dim\operatorname{Hom}(W,V\otimes u)$. Frobenius reciprocity and $u\cong\bar u$ give $F=F^{\mathsf T}$. Every irreducible occurs in a tensor power of $u$, so the graph of $F$ is connected. The vector $(1,\ldots,1)$ is fixed by $u$, whence $F_{VV}\ge1$ for every $V$. Thus $F$ is primitive, with Perron-Frobenius eigenvalue $d$ and 
Perron-Frobenius vector $(\dim V)_V$.

Let $A$ denote the adjacency matrix of the principal graph, and let $T$ be the even block of $A^2$. The principal graph is connected and $T$ has positive diagonal, so $T$ is primitive. Its Perron-Frobenius eigenvalue is $[M_H:N_H]=d$, with eigenvector $(\mu_v)_{v\ {\rm even}}$ normalized by $\mu_*=1$. Counting copies of the trivial corepresentation $\mathbf1$ in $u^{\otimes n}$ gives
\[
 (F^n)_{\mathbf1,\mathbf1}=\dim\operatorname{Fix}(u^{\otimes n})=(T^n)_{**}.
\]
For either symmetric matrix, its $n$th power divided by $d^n$ converges to the orthogonal projection onto its Perron-Frobenius vector. Since $\dim\mathbf1=\mu_*=1$, the two diagonal limits give
\begin{equation}\label{eq:gdimension}
 \frac{1}{\dim_\C A_H}=\frac{1}{\sum_V(\dim V)^2}
 =\frac{1}{\sum_{v\ {\rm even}}\mu_v^2}=\frac{1}{G}.
\end{equation}
Thus $m_H\mid h\mid G^3$, and Theorem~\ref{th:moment} puts all $Z_r(H)$ and all entries of $dQ_C$ in $\Z[\zeta_{G^3}]$. Every subfield of $\Q(\zeta_{G^3})$ is unramified outside the prime divisors of $G$ \cite[Ch.~2]{Wash}.
\end{proof}

The partial-trace identity also gives $(Q_C)_{ab}\in\Q(\operatorname{Spec}W_H)$ at arbitrary depth. Indeed, $(\operatorname{rank}E_\eta)\eta$ is a sum of eigenvalues of $W_H$. Qualitative cyclotomicity follows by another argument. Applied to the column quotients of $H^{\mathsf T}$, which are our row quotients, the construction of Jaeger, Matsumoto and Nomura identifies the $B_C$ as the adjacency matrices of the dual scheme \cite[Theorems~1 and~5(iii)]{JMN}. The Krein parameters of this scheme are integers, being intersection numbers of its dual, so Munemasa's splitting-field theorem puts its eigenvalues in a cyclotomic field \cite{Mu}. At fixed order $d$ there are only finitely many possible $Q_C$, since each entry is $1/d$ times an eigenvalue of a $d\times d$ zero--one matrix.

Taking the traces of the diagonal blocks gives $Z_1(H)=d$ and $Z_2(H)=\sum_{i,k}\Tr(P_{ik}P_{ki})$. Since $\xi_{ki}=\overline{\xi_{ik}}$, the latter identity gives \eqref{eq:zintro}.

\medskip
\noindent\textbf{An order-four example.} To illustrate Theorem~\ref{th:moment}(i) and Corollary~\ref{cor:global}, consider the family
\[
 H(q)=\begin{pmatrix}1&1&1&1\\1&-1&q&-q\\1&1&-1&-1\\1&-1&-q&q\end{pmatrix},
 \qquad |q|=1 .
\]
Burstein showed that its subfactor has finite depth exactly when $q$ is
a root of unity \cite{Bu}. To recover the necessary condition from
Theorem~\ref{th:moment}, let $\varepsilon_1,\ldots,\varepsilon_4$ be the standard basis
of $\C^4$ and put
\[
 a=(\varepsilon_1-\varepsilon_3)\otimes(\varepsilon_1-\varepsilon_2),
 \qquad
 b=(\varepsilon_2-\varepsilon_4)\otimes(\varepsilon_3-\varepsilon_4).
\]
Substituting the row-quotient projections into $W_{H(q)}=(P_{ij})$
gives $W_{H(q)}a=qb$ and $W_{H(q)}b=qa$. Thus $a+b$ is an eigenvector
with eigenvalue $q$. At finite depth $W_{H(q)}$ has finite order, so
$q$ must be a root of unity.

Squaring the entries gives two copies of each of the rows $(1,1,1,1)$ and $(1,1,q^2,q^2)$. Thus $Z_2(H(q))=8+2|1+q^2|^2$. At $q=(3+4i)/5$ this gives $272/25$, a rational number which is not an algebraic integer. The spin subfactor has infinite depth by Theorem~\ref{th:moment}(i), even though the matrix entries lie in $\Q(i)$.

If $q=\zeta_p$ for a prime $p\ge5$, the same computation gives $12+2(\zeta_p^2+\zeta_p^{-2})$. This generates the maximal real subfield of $\Q(\zeta_p)$, in which $p$ ramifies. Since $q=\zeta_p$ gives finite depth, $p$ must divide the global index, although the Jones index is always four.

\begin{proposition}[Conference matrices]\label{prop:conference}
Let $C$ be a real symmetric conference matrix of order $d>2$; that is, $C$ has zero diagonal, off-diagonal entries in $\{1,-1\}$, and $C^2=(d-1)I$. Then the complex Hadamard matrix $H=I+iC$ has an infinite-depth spin subfactor.
\end{proposition}
\begin{proof}
We have $HH^*=I+C^2=dI$. Multiplying rows or columns of $H$ by phases does not change any $P_{ij}$, so we may apply Theorem~\ref{th:moment} after dephasing. Write $S=(H_{ia}^2)_{i,a}=2I-J$, where $J$ is the all-ones matrix. By \eqref{eq:zintro},
\[
 Z_2(H)=d^{-2}\Tr\bigl((SS^*)^2\bigr)
       =\frac{(d-2)^4+16(d-1)}{d^2}
       =d^2-8d+24-\frac{16}{d}.
\]
Orthogonality of two rows of $C$ implies that $d$ is even. After simultaneous sign changes of rows and columns, the off-diagonal entries of the first row and column are all one. The remaining symmetric core has row sums zero, so each row has $(d-2)/2$ positive entries. Their total number is even by symmetry; since the core has $d-1$ rows, we obtain $d\equiv2\pmod4$. Thus $d>2$ does not divide $16$, and $Z_2(H)$ is rational but not integral. This contradicts Theorem~\ref{th:moment}(i) at finite depth.
\end{proof}

\begin{corollary}\label{cor:conductor}
Fix integers $d\ge2$ and $D\ge1$. There are only finitely many dephased complex Hadamard matrices $H$ of order $d$ for which the principal graph of $N_H=A_{\infty,0}\subset M_H=A_{\infty,1}$ has depth at most $D$. There is an integer $m=m(d,D)\ge1$ such that every $W_H$ in this family has order dividing $m$. All its power traces and all entries of $dQ_C$ then belong to $\Z[\zeta_m]$.
\end{corollary}
\begin{proof}
The order $d$ fixes the four corners and their inclusion matrices. By Theorem~\ref{th:main} only finitely many isomorphism classes occur. Each gives at most $(d!)^2$ matrices $H$, since ordering the two sets of minimal projections determines $H_{ij}=d^3\tau(p_iq_jp_1q_1)$. Take $m$ to be the least common multiple of their orders $\operatorname{ord}(W_H)$, or $1$ for the empty family, and apply Theorem~\ref{th:moment}.
\end{proof}


\end{document}